\pdfoutput=1
\documentclass{article}
\usepackage[margin=25mm]{geometry}
\usepackage{amsmath}
\usepackage{amsfonts}
\usepackage{amssymb}
\usepackage{graphicx}
\usepackage[utf8]{inputenc}
\usepackage[T2A]{fontenc}
\usepackage[ukrainian,russian,english]{babel}
\usepackage[unicode, pdftex]{hyperref}
\usepackage{xcolor}
\usepackage[all]{xy}

\newtheorem{theorem}{Theorem}

\title{ Gradient vector fields of codimension one on the 2-sphere with at most ten singular points }

\author{Svitlana Bilun, Bohdana Hladysh, 
Alexandr Prishlyak and Vladislav Sinitsyn}
\begin{document}

\maketitle
\begin{abstract}
   We describe all possible (469) topological structures of сodimension one gradient vector fields  on the 2-sphere with at most ten singular points. To describe structures, we use a graph whose edges are one-dimensional stable manifolds. The saddle-node singularity is specified by selecting a pair of vertices-edge or edge-face, and the saddle connection is specified by a T-vertex.
\end{abstract}

\section{Introduction}

In this paper, we consider gradient vector fields on a sphere. Since the function increases along each trajectory, the field has no cycles and polycycles. In general position, a typical gradient field (a vector field of the codimension 0) is a Morse field (Morse-Smale field without closed trajectories). In typical one-parameter families of gradient vector fields, two types of bifurcations are possible: a saddle-node and a saddle connection. The corresponding vector fields at the time of the bifurcation are fields of codimension 1. In our case, they completely determine the topological type of the bifurcation.
To classify Morse fields, a cell complex structure (diagram) is often used, in which cells of dimension n are stable manifolds of singular points with Morse index equal to n. We apply this approach to the classification of vector fields of codimension 1 on the 2-sphere. The cell structure on the 2-sphere is determinated by the imbedded graph, which is the 1-skeleton.  

Without loss of generality, we assume that the number of singular points does not increase under bifurcation (as the parameter increases). The saddle-node bifurcation is defined by a pair of cells corresponding to the singular points participating in the bifurcation. We mark this pair on the diagram with a green arrow or a green triangle. A saddle-node bifurcation in the diagram corresponds to a point of degree 3 (T-vertex), where two edges (half-edges) are opposite and the third is perpendicular to them.

Then, the separatrix that connects the saddle with the node (source or sink) contracts to a point under the saddle-node bifurcation.

Selection of pairs of cells of neighboring dimensions is used in the discrete Morse theory \cite{for98}.

Topological invariants of functions were constructed in the papers of Kronrod \cite{Kronrod1950} and Reeb \cite{Reeb1946} for oriented maniofolds, in \cite{lychak2009morse} for  non-orientable two-dimensional manifolds and in   \cite{Bolsinov2004, hladysh2017topology, hladysh2019simple} for manifolds with boundary, in \cite{prishlyak2002morse} for non-compact manifolds. 

In general, Morse vector fields (Morse-Smale vector fields without closed orbits) are gradient field of Morse functions. If we fix the value of functions in singular points the field determinate the topological structure of the function \cite{lychak2009morse, Smale1961}. Therefore, Morse--Smale vector fields classification is closely related to the classification of the functions.

Topological classification of smooth function on closed 2-manifolds was also investigated in \cite{bilun2023morseRP2,  hladysh2019simple, hladysh2017topology,  prishlyak2002morse, prishlyak2000conjugacy,  prishlyak2007classification, lychak2009morse, prishlyak2002ms, prish2015top, prish1998sopr,  bilun2002closed,  Sharko1993, bilun2013def}, on 2-manifolds with the boundary in \cite{hladysh2016functions, hladysh2019simple, hladysh2020deformations} and on closed 3-manifolds in  \cite{prishlyak1999equivalence}.

In \cite{Kybalko2018, Oshemkov1998, Peixoto1973, prishlyak1997graphs, prishlyak2020three, akchurin2022three, prishlyak2022topological, prishlyak2017morse,  kkp2013,  prish2002vek,  prishlyak2021flows,  prishlyak2020topology,   prishlyak2019optimal, prishlyak2022Boy}, 
the classifications of flows on closed 2- manifolds and 
\cite{bilun2023discrete, loseva2016topology, prishlyak2017morse, prishlyak2022topological, prishlyak2003sum, prishlyak2003topological, prishlyak2019optimal} on manifolds with the boundary were obtained.
Topological properties of Morse-Smale vector fields on 3-manifolds was investigated in \cite{prish1998vek,  prish2001top, Prishlyak2002beh2, prishlyak2002ms,   prish2002vek, prishlyak2005complete, prishlyak2007complete, hatamian2020heegaard, bilun2022morse, bilun2022visualization}.
%prishlyak2003regular,

The topological  invariants of graphs and their embeddings in the 2-manifolds can be founded in \cite{prishlyak1997graphs, Harary69, HW68, GT87}.

%The topological properties of the projective plane are described in many geometry textbooks, but we would also recommend \cite{Bilun22Projective}.

In view of the bijection between vector fields and flows on closed manifolds, we identify these two concepts.

The purpose of this paper is to describe all possible topological structures of the gradient flow structure of dimension 1 on a sphere with no more than four saddles (a saddle-node point is also considered a saddle).

%The first section gives the basic definitions and topological properties of discrete Morse functions and discrete gradient vector fields.

%In the second section we describe all possible structures of optimal discrete vector flows on a two-dimensional disk, in the third section -- on a two-dimensional sphere, in the fourth section -- on the cylinder, and in the fifth section -- on a Mobius strip.

%\subsection{Постановка задачі} 
%Метою даної роботи є описати всі можливі топологічні структури структура градієнтних потоків короз\-мір\-но\-сті 1 на сфері з не більш ніж чотирма сідлами (точка сідло-вузол також вважається сідлом). 
 
%\subsection{Основні поняття}
%Сепаратрисна діаграма.

\section{Typical one-parameter bifurcations of gradient flows on a sphere}

Typical vector fields on compact 2-manifolds are Morse-Smale fields. Among the gradient fields are Morse fields or Morse-Smale gradient-like fields that do not contain closed trajectories. They satisfy three properties:

1) singular points are nondegenerate;

2) there are no saddle connections;

3) $\alpha$-limit ($\omega$-limit) set of each trajectory is a singular point.

In typical one-parameter field families, one of these conditions is violated. Violation of the first condition leads to a saddle-node bifurcation, and the second to the appearance of a saddle connection. The third condition cannot be violated, because according to the Poincaré-Bendixon theorem, the $\alpha$-limit ($\omega$-limit) set of every trajectory on the sphere is either a singular point, or a cycle, or a polycycle. Since gradient fields have no cycles and polycycles, this set is a singular point.

\subsection{Saddle-node}
The saddle-node bifurcationis shown in Fig. 1 if the node is the source.
\begin{figure}[ht!]
\center{\includegraphics[width=0.95\linewidth]{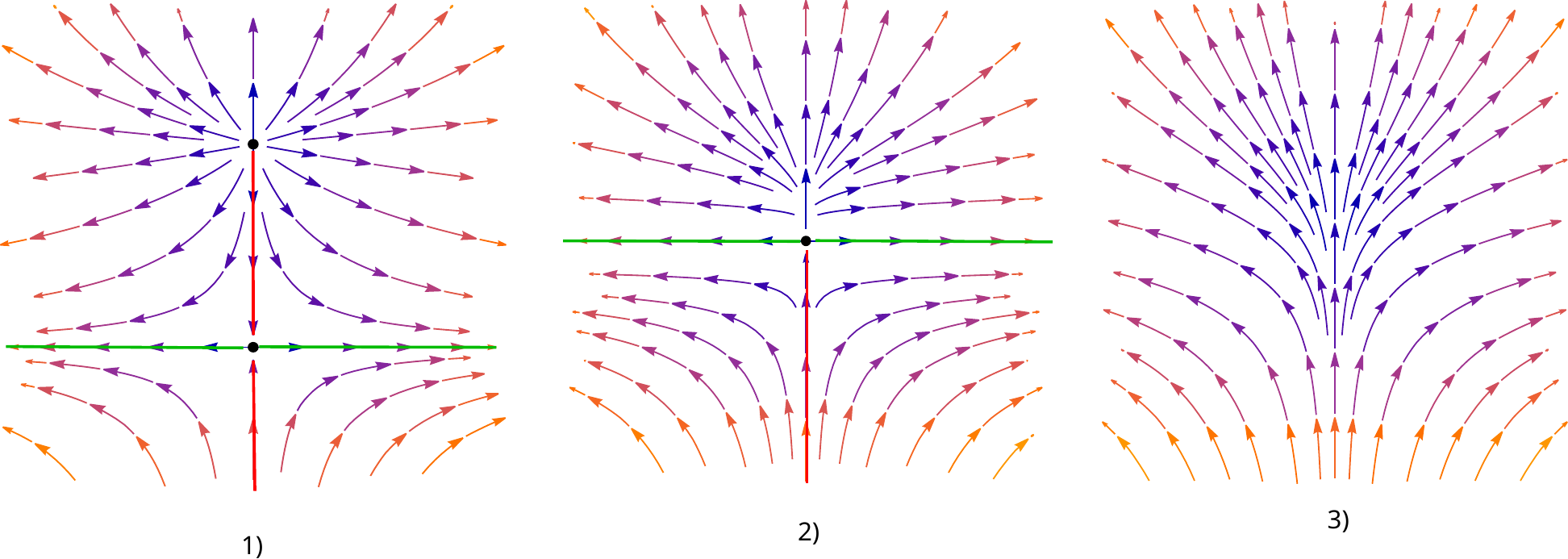}
}
\label{bifsn}
\caption{saddle-nod bifurcation
}
\end{figure}
It can be described by the equation $V(x,y,a)=\{x,y^2+a\}$ if the node is the source and the equation $V(x,y,a)=\{-x,- y^2a\}$ if the node is a sink. Here, $a$ is a parameter. If $a<0$ we get the flow before the bifurcation, if $a>0$ we get the flow after the bifurcation, and if $a=0$ we get the flow at the moment of the bifurcation (the flow of codimensionality 1). The saddle-sink bifurcation can be obtained from the saddle-source rotation of the direction on all trajectories, that is, by taking the reversed flow. Note that the saddle-node bifurcation corresponds to a typical one-parameter-method catastrophe (deformation) of functions: $f(x,y,a)=x^2+y^3+ay$.

Consider a stable manifold (red separatrix in the figures), for a singular saddle-source point. It cannot start at this particular point because the gradient field contains no loops. This means that a stable saddle manifold participating in a saddle-source bifurcation does not form a loop. Therefore, in order to determine the saddle-source bifurcation on the cellular partition given by the stable multiples of the Morse field before the bifurcation, it is necessary to select a 1-cell that is not a loop and one of its ends. In the following figures, we will do this with the help of a blue arrow or blue triangles at the ends of the ribs.

Similarly to the above, a saddle-sink bifurcation is possible if the unstable manifold of the saddle before the bifurcation does not form a loop. This means that the corresponding 1-cell is contained in the boundary of two different 2-cells. To set such a bifurcation, you need to select a pair: 1 cell and one of the two adjacent 2-cells. In the pictures, we will do this with the help of a blue arrow (or a blue triangle) directed from the 1-cell towards the corresponding 2-cell.

Note that the bifurcation of the saddle flow in Morse's theory corresponds to the operation of reducing critical points or reducing complementary handles.

\subsection{Saddle connection}
The bifurcation of the saddle connection is shown in Fig. 2. It can be described by the equation $V(x,y,a)=\{x^2-y^2-1,-2xy+a\}$.
\begin{figure}[ht!]
\center{\includegraphics[width=0.95\linewidth]{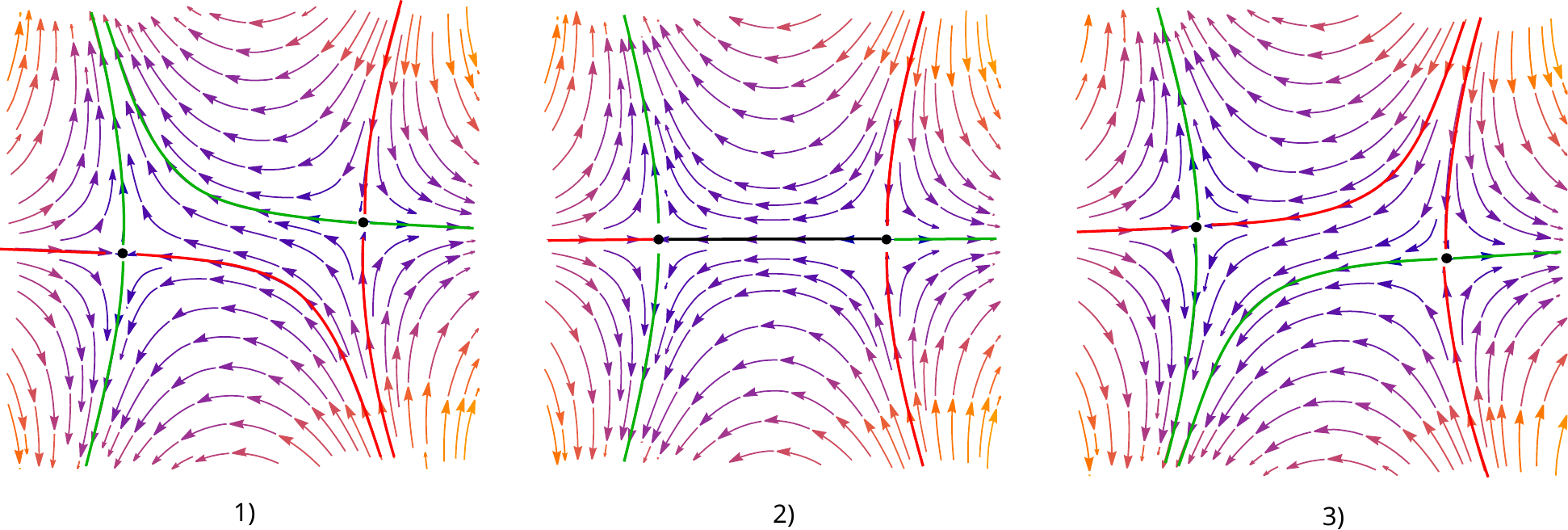}
}
\label{bifsc}
\caption{saddle connection
}
\end{figure}
We describe one of the possible situations when it appears in typical one-parameter families of gradient fields. Let $p,q$ be saddle critical points of the function $f$, $f(p)<f(q)$, $[\frac{f(p)+f(q)}{2}-\varepsilon, \frac{f(p)+f(q)}{2}+\varepsilon]$ does not contain critical values of $f$, $L$ of the line component of level $f^(-1)(\frac{f(p) +f(q)}{2})$, $u$ is the separatrix of $p$, $v$ is the separatrix of $q$, which intersect $L$.
Then according to Morse theory, $f^{-1}([\frac{f(p)+f(q)}{2}-\varepsilon, \frac{f(p)+f(q)}{2} +\varepsilon])$ is homeomorphic to the cylinder $S^1 \times [0,1]$ and in the coordinates $(s,t)$ ($s$--polar angle on the circle) on it, the trajectories are given by segments $s= \text{const}$. Consider the twisting of a cylinder given by the formula $(s,t) \to (s+a,t)$ ($a$-- deformation parameter). At the same time, the trajectories  turn into helical lines. With a continuous change of the parameter $a$, the point of intersection of the trajectory $u$ with the circle $frac{f(p)+f(q)}{2}+\varepsilon$ continuously rotate along it. The intersection point of the trajectory $v$ with this circle remain unchanged. According to the intermediate value theorem, there is a moment in time $a$ at which these points coincide.Then the bifurcation of the saddle connection occur  at this point of time $a$.
Note that the construction described above can be implemented by changing the Riemannian metric (changing angles considered perpendicular), without changing the function itself.

In the figuress, the bifurcation of the saddle connection will be indicated by a rib that is perpendicularly glued to the middle of the other rib (in the form of the letter T). A glued edge (1-cell) containing a separatrix connection  be called perpendicular. Therefore, a graph containing a T-vertex and a perpendicular edge corresponds to a flow of codimensionality 1.

In Morse's theory, bifurcation of the saddle connection corresponds to the operation of adding handles (sliding one handle over another handle of the same index).

\subsection{Distinguishing graph }

Let's construct a graph embedded in a sphere, the vertices of which are sources, saddle-node points, and T-vertices, and the edges are stable one-dimensional manifolds. On the graph, we select a T-vertex, if it is present and has a perpendicular edge in it. And in the case of a saddle-source, we  select the corresponding vertex-edge pair, and in the case of a saddle-sink, we select an edge-face pair. A spherical graph with a selected element (T-vertex and edge, vertex-edge or edge-face pair) is called a distinguishing graph of a vector field or simply a field graph. Two such graphs are equivalent if there is an isomorphism of them as spherical graphs that preserves the selected element.

\begin{theorem}
Two vector fields of codimensionality 1 on the sphere are topologically equivalent if and only if their distinguishing graphs are equivalent. Each spherical graph with a selected pair is a distinguishing graph of a vector field of codimensional 1.
\end{theorem}
\textbf{Proof.} The necessity follows directly from the construction of the distinguishing graph. To prove sufficiency, let's construct a separatrix diagram for the graph. To do this, consider a dual graph. If there is no T-vertex, then the selected pair defines a half-edge on the initial or dual graph. Let's squeeze it into a point. The remaining half-edges are flow separatrixes. Their direction is determined from the vertices of the initial graph to the vertices of the dual graph. So, we got a separatrix diagram of a vector field, which specifies its topological type. In the case of a T-vertex, for the perpendicular edge, we still construct a double edge, and for the other two edges incident to the T-vertex, we construct an edge that is a continuation of the perpendicular edge. Therefore, the T-vertex turn into a vertex of valence 4 (saddle). Further considerations are the same as in the case without a T-vertex.

Since the separatrix diagram of a graph is defined uniquely, such a graph defines a vector field.

\section{Fields with one and two saddles. }

Let us describe the structure of all possible vector fields of codimensionality 1 with one and two saddles. For saddle-node bifurcations, we define the flows to the bifurcations by their graphs. The number of ribs is equal to the number of saddles. There are two graphs with one edge: 1) $G^1_1$ -- with two vertices (segment -- Fig. 3 1)), 2) $G^1_2$ -- with one vertex (loop --  Fig.3. 7) ).
\begin{figure}[ht!]
\center{\includegraphics[width=0.6\linewidth]{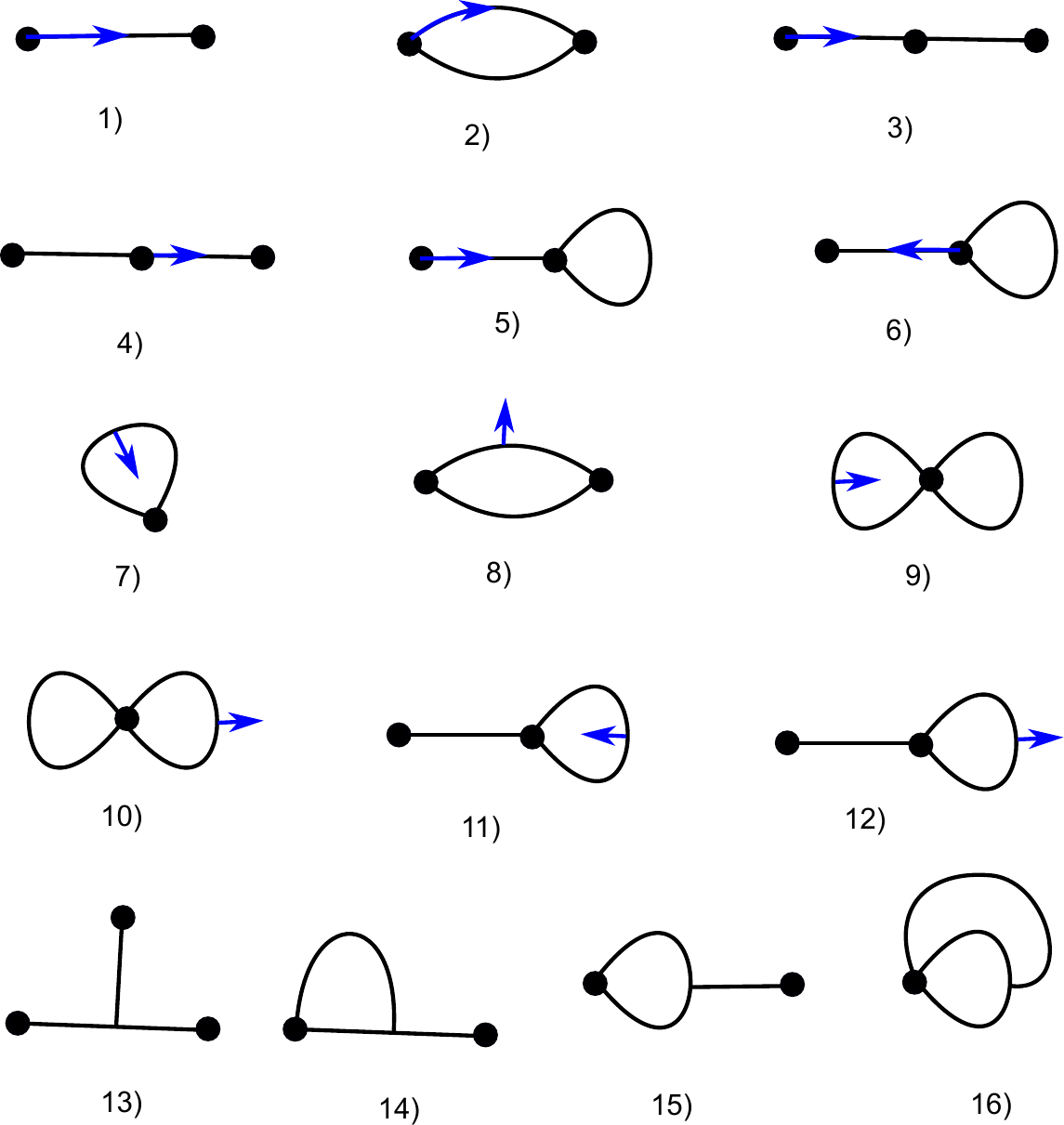}
}
\label{bif12}
{
\caption{Flows of codimensionality 1 with one and two saddles. 1)-6) saddle-source, 7)-12) saddle-sink, 13)-16) saddle connection. Graphs 1) $G^1_1$, 2) $G^2_1$, 3) $G^2_2$, 5) $G^2_3$, 7) $G^1_2$ 9) $G^2_4$. }
}
\end{figure}

Only four graphs are possible with two edges: 1) $G^2_1$ -- a pair of multiple edges, 2) $G^2_2$ -- a chain of length 2, 3) $G^2_3$ -- a segment with a loop, 4) $ G^2_4$ -- a pair of loops.

On $G^1_1$ there is a single (with accuracy to homeomorphism) saddle-source bifurcation (see Fig.3. 2). Saddle-sink  is not possible because $G^1_1$ does not split the plane. On $G^1_2$, on the contrary, there is a saddle-sink bifurcation, and a saddle-source bifurcation does not exist because the only edge is a loop.
Note that $G^1_2$ is dual to $G^1_1$, the corresponding vector fields are obtained from each other by multiplication by $-1$. When the direction of the vector field changes, the source turns into a sink, so the saddle-source bifurcation turns into a saddle-sink.
Bifurcation of the saddle connection for flows with one saddle is not possible, since two saddles take part in it.

On $G^2_1$, all possible saddle-source pairs are symmetric, so there is only one such bifurcation 2).
Given the symmetry of $G^2_2$, there are two saddle-node bifurcations on it: 3) -- with an extreme vertex, 4) -- with a central vertex. On the graph $G^2_3$, the loop cannot participate in the saddle-source bifurcation. For another segment edge, there are two options for choosing its end. So, we have two saddle-source bifurcations 8) and 9). Saddle-source bifurcation is not possible on a pair of loops. If we consider the reverse flows, then the resulting 5 saddle-source bifurcations will turn into saddle-sink bifurcations on the dual graphs shown in Fig.3. 8)-12).

If both saddles participate in the bifurcation of the saddle connection, then one of them will be normal, and the other - perpendicular. If the normal edge is a segment, then there are two options for the perpendicular edge: 13) and 14). If the regular edge is a loop, then we also have two options 15) and 16).

We went through all the possible options, and therefore it is fair

\begin{theorem}
The sphere has the following structures of codimensionality 1 flows:
\begin{itemize}
\item
two flow structures with three singular points (1 and 7);
\item
there are no streams with four points;
\item
ten structures with five points (2--6 and 8--12);
\item
four structures with six points (13--16).
\end{itemize}
\end{theorem}

\section{Fields with three saddles}
In the following, we will describe only saddle-source bifurcations, since they specify all possible saddle-sink bifurcations on dual graphs. Let us describe all possible graphs with three edges, see Fig.4.

  \begin{figure}[ht!]
\center{\includegraphics[width=0.8\linewidth]{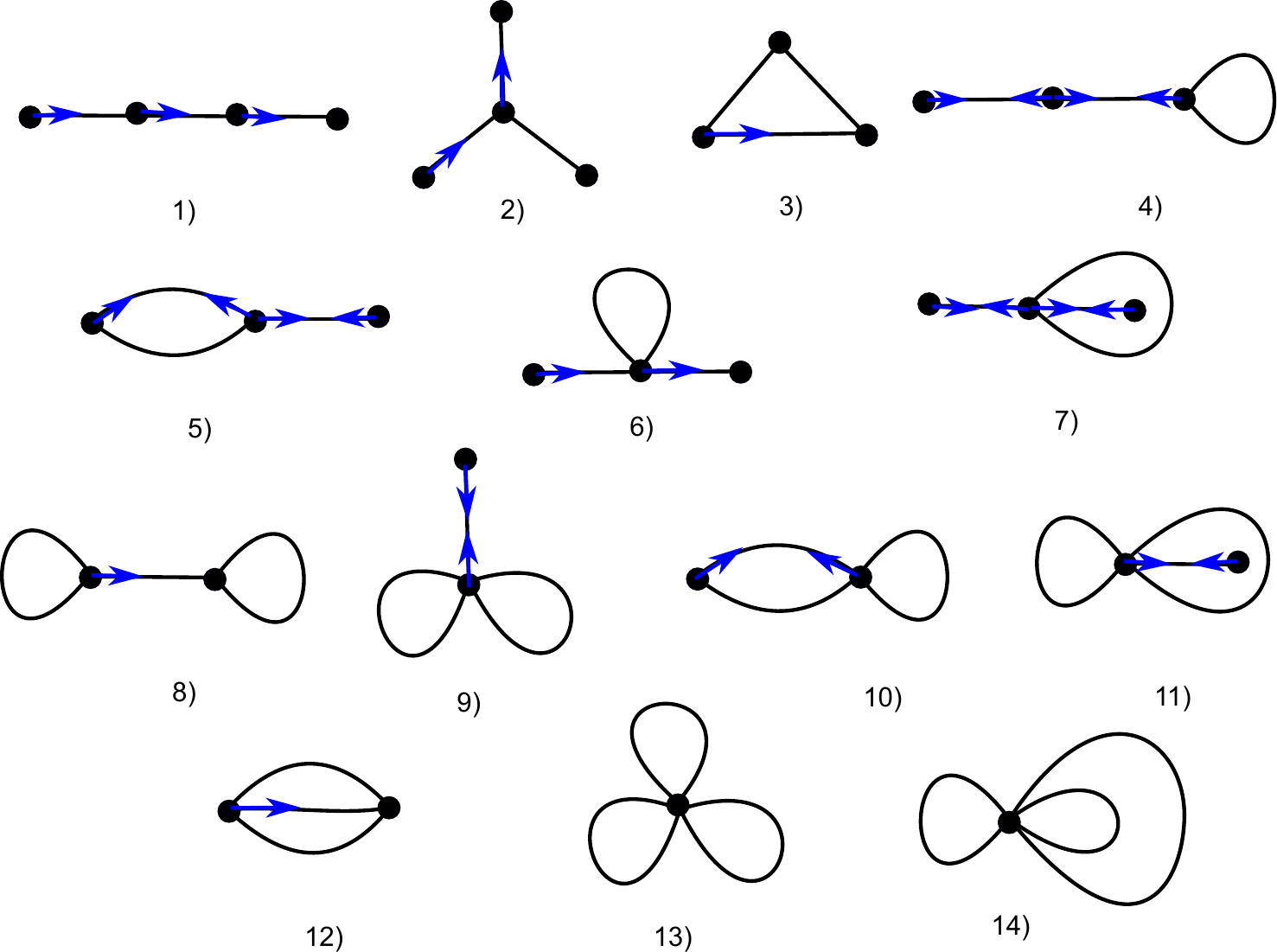}
}
\label{bif3a}
\caption{Fields of codimensionality 1 with three saddles: the saddle-source case
}
\end{figure}
  \begin{figure}[ht!]
\center{\includegraphics[width=0.90\linewidth]{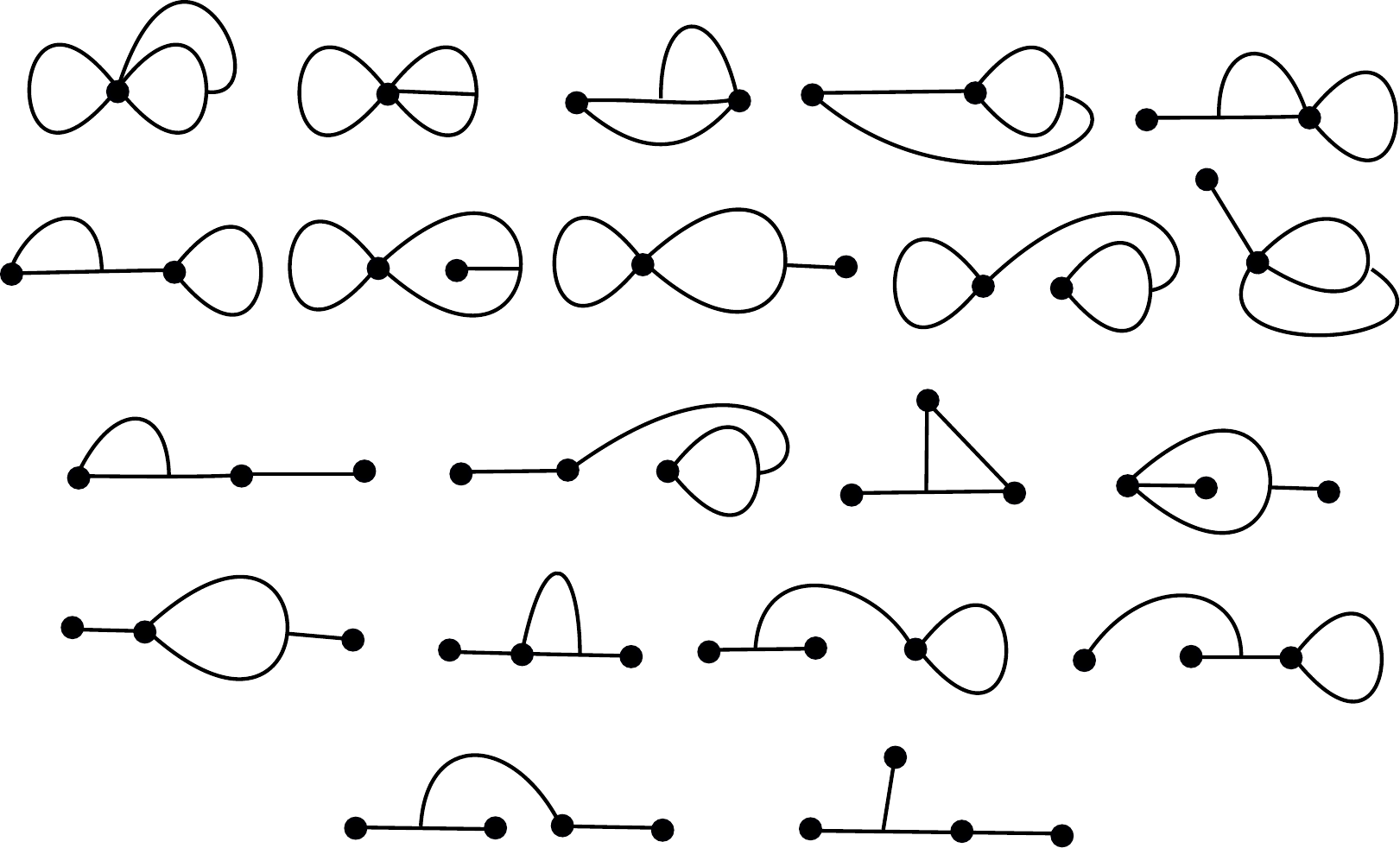}
}
\label{bif3b}
\caption{Flows of codimensionality 1 with three saddles: a peculiarity of the saddle connection
}
\end{figure}
Let a graph contain four vertices, then it is a tree.
There can be only two such trees: the chain $G^3_1$ and the Y-graph $G^3_1$. If the graph has three vertices, then they are connected by a chain. Let's number them so that the second vertex is the middle one on this chain. Let us first consider the case when the third edge is not incident to the second vertex. If it connects the first and third vertices, then we get a triangle $G^3_3$, and if it has one vertex at its ends, then the graph $G^3_4$ is two chains with a loop at the end. If the third edge starts from vertex 2, then we have two options: 1) it ends at the extreme vertex -- graph $G^3_5$ or at vertex 2 -- graphs $G^3_6$ and $G^3_7$.
Now let's write all binary graphs in reverse order. We will obtain graphs $G^3_8$ -- $G^3_12$ with two vertices and graphs $G^3_{13}$ and $G^3_{14}$ with one vertex.

In Fig.5, each of the graphs $G^3_k$ is depicted above its number $k$).

All possible saddle-source bifurcations of these graphs are depicted with blue arrows.

In Fig.5, all possible bifurcations of the saddle connection are shown
with three saddles.

So, fair

\begin{theorem}
The sphere has the following structures of codimensionality 1 flows:
\begin{itemize}
\item
56 structures with seven singular points;
\item
20 structures with eight singular points.
\end{itemize}
\end{theorem}

\section{Flows with four saddles}

\begin{figure}[ht!]
\center{\includegraphics[width=0.8\linewidth]{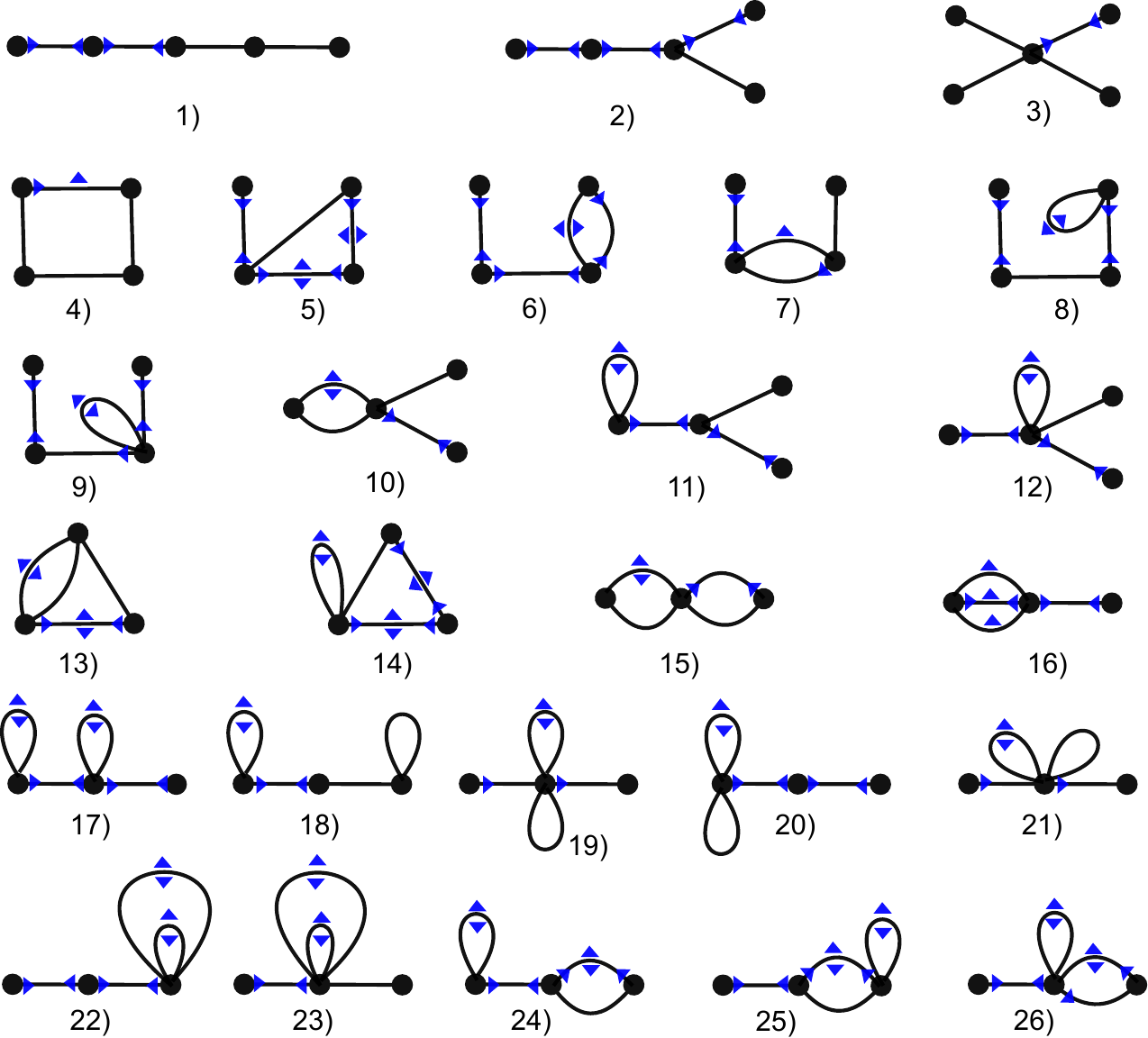}
}
\label{bif4sn}
\caption{Vector fields with four saddles: the saddle-node case
}
\end{figure}

Graphs with four edges (see Fig.\ref{bif4sn}) are of the following types:

a) with five vertices (1--3),

b) with four vertices (4--12),

c) with three vertices (13--26),

d) with two vertices (dual to 4--12),

e) with one vertex (dual to 1--3).

Each graph of type e) is dual to a graph of type a). Graphs of type d) are dual to graphs of type b). Therefore,  only graphs with five, four or three vertices and all their possible saddle-node bifurcations using blue triangles are shown in Fig.6.

Saddle-node bifurcations give inverse saddle-node bifurcations on a dual graph. Therefore, in order to find the total number of flows of codimensional 1 of this type, we multiply the number of fields (64) given by graphs 1--12 by 2 and add the number of fields (89) given by graphs 13--26. Finally we have:

\begin{theorem}
On the sphere there are 217 different codimensional 1 vector field structures with nine singular points.
\end{theorem}

Next, we consider bifurcations of the saddle connection in flows with four saddles (10 singular points). Two cases are possible for an edge perpendicular to a T-vertex: 1) it does not break the graph, 2) it breaks the graph.

In the first case, if it is removed, we obtain a connected graph with three edges (if two edges that are incident T-vertices are connected into one edge). The removed edge is defined by a perpendicular vector in the T-vertex and one of the corners of the face in which it is directed. We add one to these numbers, which corresponds to the case when the perpendicular edge has a vertex of degree 1 at one of its ends. All possible such cases (130) are shown in Fig.7.

  \begin{figure}[ht!]
\center{\includegraphics[width=0.75\linewidth]{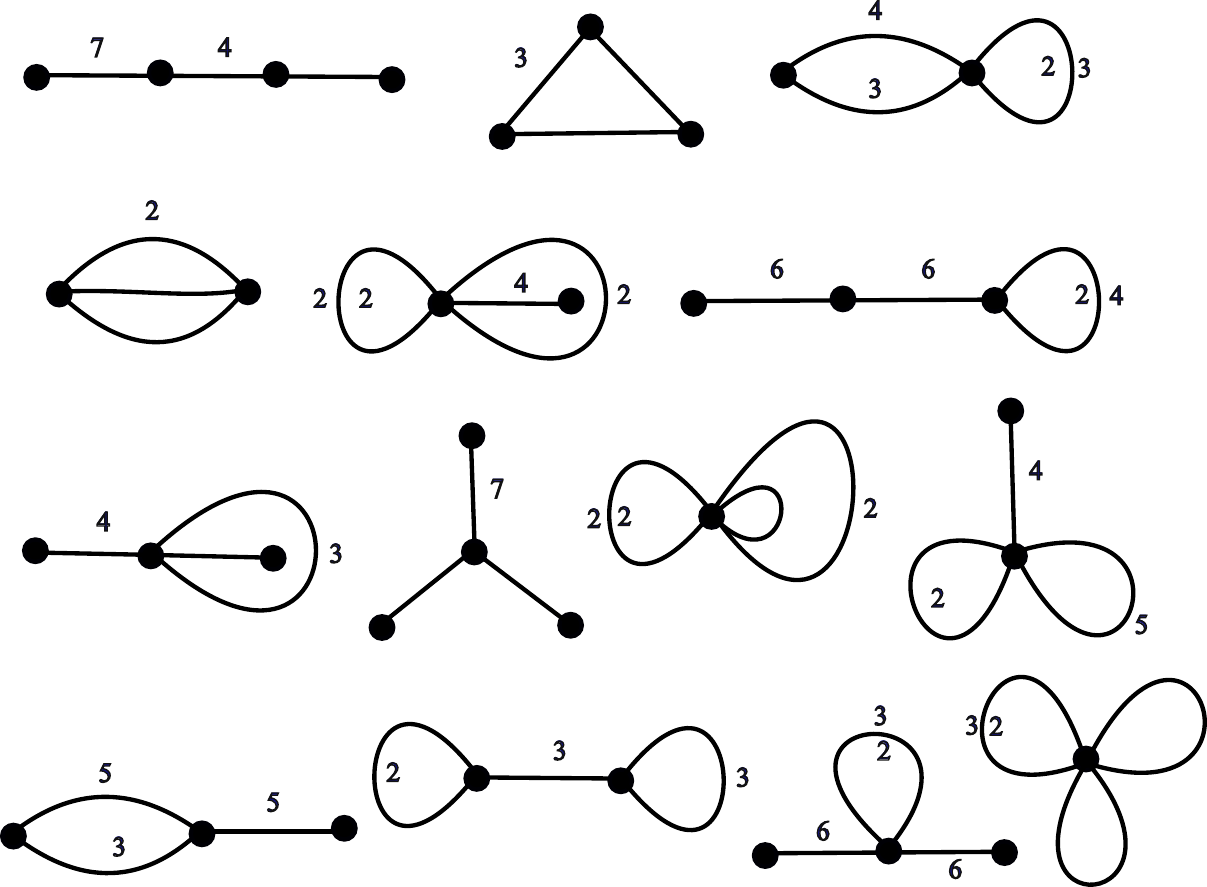}
}
\label{bif4sc}
\caption{Vector fields with four saddles: a case of the saddle connection
}
\end{figure}

Consider the case when a perpendicular edge divides the graph into two subgraphs. Let's glue two edges at the vertex T into one and denote this subgraph by $K$, and the other subgraph by $L$. Then  $K$ can contain one, two or three edges. In the latter case, the graph $L$ has no edges and only one vertex, which we considered earlier. If the graph $K$ has one edge ($G^1_1$ or $G^1_2$), then in each of them the perpendicular vector is determined uniquely (up to homeomorphism), that is, we have 2 options. On graphs with two edges, there are 8 different possibilities to choose the corner at the vertices ($G^2_1$ -- 1, $G^2_2$ -- 2, $G^2_3$ -- 3, $G^2_4$ -- 2) . Therefore, $2 \times 8 = 16$ fields of this type are possible.

If $K$ has two edges, then there are 7 ways for the perpendicular vector to the edges: $G^2_1$ -- 1, $G^2_2$ -- 1, $G^2_3$ -- 3, $G^2_4$ -- 2. The graph $L$ contains one edge and in each of two such graphs the angle is specified uniquely. Therefore, we have $7 \times 2 = 14$ fields of this type.

So, fair
\begin{theorem}
There are 160 different flow structures of codimensionality 1 with ten singular points on the sphere.
\end{theorem}

\section*{Conclusion}

All possible structures of flows of codimensionality 1 (typical bifurcations) on a two-dimensional sphere with no more than four saddles (10 singular points) were found: with three singular points -- 2 structures, with four -- none, with five points -- 10, with six -- 4, with seven -- 56, with eight -- 20, with nine -- 217 and with ten points -- 160 structures.
We hope that the research carried out in this paper can be extended to other surfaces: spheres with holes, surfaces of a larger genus, in particular non-orientable ones.

%\bibliographystyle{plain}
%\bibliography{prish}

\textsc{Taras Shevchenko National University of Kyiv}

Svitlana Bilun  \ \ \ \ \ \ \ \ \ \ \ \
\textit{Email address:} \text{ bilun@knu.ua}   \ \ \ \ \ \ \ \ \ \
\textit{ Orcid ID:}  \text{0000-0003-2925-5392}

Bohdana Hladysh \ \ \ \ \ \ \ 
\textit{Email:} \text{ bohdanahladysh@gmail.com,}  \ 
\textit{ Orcid ID:} \text{0000-0001-8935-1453}

Alexandr Prishlyak \ \ \ \ \ 
\textit{Email address:} \text{ prishlyak@knu.ua} \ \ \ \  \
\textit{ Orcid ID:} \text{0000-0002-7164-807X}

Vladislav Sinitsyn \ \  \ \ \ \ \  \textit{Email:} \text{vlad.sinitsyn.30@gmail.com }    \ \
\textit{ Orcid ID:} \text{0009-0000-0928-0676}

\end{document}